\documentclass{article}
\usepackage{amsmath, amssymb, amsthm, enumerate}
\title{Some implications of the Gessel identity}
\author{Claire Levaillant}

\newcommand{\sa}{\sum_{a=1}^{p-1}}
\newcommand{\om}{\omega}

\newcommand{\mpt}{\;\text{mod}\,p^3}
\newcommand{\mpd}{\;\text{mod}\,p^2}
\newcommand{\mpu}{\;\text{mod}\,p}

\newcommand{\mc}{\mathcal{C}}

\newcommand{\mb}{\mathcal{B}}

\newcommand{\md}{\mathcal{D}}
\newcommand{\mh}{\mathcal{H}}

\newcommand{\mt}{\mathcal{T}}

\newcommand{\as}{A^{\star}}

\newcommand{\mpq}{\;\;\text{mod}\,p^4}

\newcommand{\mbt}{\mb\mb\mb}

\newcommand{\Om}{\Omega}

\begin{document}
\maketitle
\begin{center}Abstract.\end{center}

\noindent We generalize the congruences of Friedmann-Tamarkine ($1909$), Lehmer $(1938)$, Ernvall-Mets\"ankyla $(1991)$ on the sums of powers of integers weighted by powers of the Fermat quotients to the next Fermat quotient power, namely to the third power of the Fermat quotient. Using this result and the Gessel identity ($2005$) combined with our past work $(2021)$, we are able to relate residues of some truncated convolutions of Bernoulli numbers with some Ernvall-Mets\"ankyla residues to residues of some full convolutions of the same kind. We also establish some congruences concerning other related weighted sums of powers of integers when these sums are weighted by some analogs of the Teichm\"uller characters.

\section{Introduction}

\subsection{Summary}

\indent

In $1991$, Reijo Ernvall and Tauno Mets\"ankyla generalized the K\"ummer congruences by relating modulo $p^2$ the difference $\mb_{p-1+t}-\mb_t$ 
to a sum of $t$-th powers of the first $(p-1)$ integers weighted by the squared Fermat quotients. They did so for the range $4\leq t\leq p-3$ when $t$ is even.
In their honor, we will call the first non-zero $p$-residue of $\mb_{p-1+t}-\mb_t$ the Ernvall-Mets\"ankyla residue. In this paper, we set a prime $p\geq 11$ and we restrict further the range for the even $t$. We namely impose $6\leq t\leq p-5$. We let $n$ be the integer so that $t=p-1-2n$. Hence $4\leq 2n\leq p-7$. 
The main result of the paper offers a congruence modulo $p$ relating a truncated convolution of order $2(p-1)-2n$ of Bernoulli numbers with Ervall-Mets\"ankyla residues to the full convolution of order $p-1-2n$. The latter congruence also involves a full convolution where the ordinary Bernoulli numbers are replaced with the divided Bernoulli numbers. This result can be viewed as a generalization of Theorem $1$ point $(i)$ of \cite{LEV3} which offers a similar result when dealing with convolutions of divided Bernoulli numbers. One of the three proofs presented in \cite{LEV3} uses the famous Miki identity which relates convolutions of divided Bernoulli numbers to binomial convolutions of divided Bernoulli numbers. The current proof is based on Ira Gessel's identity, namely a generalization of Hiroo Miki's identity to cubic convolutions, and on our past work on the multiple harmonic sums $\mathcal{H}_{\lbrace s\rbrace^{2n}=1;p-1}$ modulo $p^4$, amongst other things. Another important tool consists of generalizing congruences proven throughout the twentieth century and which are due to Friedmann-Tamarkine \cite{FT}, Lehmer \cite{LEM}, Miki \cite{MI} and Ernvall-Mets\"ankyla \cite{EM} in chronological order, on sums that are discussed right below and which involve the Fermat quotients in base $a$ with $1\leq a\leq p-1$. The latter generalization is interesting in its own sake. Though many mathematicians have been studying congruences concerning Bernoulli numbers and the Fermat quotients, their results tend to be sparse in time, like testified by the (non exhaustive) list below.\\
\indent Since the beginning of the twentieth century, mathematicians have been working with sums of some fixed power of the first $(p-1)$ integers involving the Fermat quotients or the squared Fermat quotients as their weights. However, sums of some fixed power of integers have been studied since way earlier. In fact, the pioneer for the study of such sums is Johann Faulhaber of Ulm $(1580-1635)$ whose computations by hand in the early $17$th century \cite{JF} drew the attention and the admiration of many mathematicians in the subsequent centuries. Two centuries after Faulhaber, in $1834$, Carl Jacobi was first to provide a rigorous proof for Faulhaber's guessed formula \cite{JA} which later became known as "the Faulhaber formula", see \cite{KN} for an excellent expository on the topic. \\
Weighted sums of powers of integers happen to be linked to the Bernoulli numbers. Conversely, the residues of the divided Bernoulli numbers, to the exception of some of them, can be expressed as residues of some weighted sums of powers of integers. Indeed, if $i$ is a positive even integer that is prime to both $p$ and $p-1$, then
$$\mb_i=\sa\Bigg(\sum_{b=1}^{p-1}\left[\begin{array}{l}\frac{b^{-1}a}{p}\end{array}\right]\frac{1}{b^{i-1}}\Bigg)a^{i-1}\mpu$$
This congruence namely follows from summing the Voronoi congruences from $1889$ \cite{VOR} and using the fact that $p$ divides a sum of powers of the first $(p-1)$ integers when these powers are not multiples of $(p-1)$.\\
Three centuries after Faulhaber, the earlier work on the sums of some fixed power of the first $(p-1)$ integers weighted by the Fermat quotients goes back to Friedmann and Tamarkine \cite{FT} in the early twentieth century. Their congruence modulo $p$ opens another path towards the study of irregular primes since it provides a criterion for irregular pairs. The search for irregular primes had been of importance on the route to the proof of Fermat's Last Theorem ever since Ernst K\"ummer had shown that Fermat's Last Theorem holds when the exponent is a regular prime \cite{KUM}, thus leaving irregular primes to be analyzed individually. While Friedmann and Tamarkine only deal with even powers, Emma Lehmer in $1938$ generalizes the Friedmann-Tamarkine congruence to odd powers and pushes the study one $p$ power further. Her congruences namely hold modulo $p^2$. It is not until $1978$ that sums of some fixed power of the first $(p-1)$ integers weighted by the squared Fermat quotients appear in the pioneering work of Hiroo Miki on his way to finding an identity relating a convolution of divided Bernoulli numbers to a binomial convolution of divided Bernoulli numbers, an identity which also involves harmonic numbers. Miki relates these sums modulo $p$ for a restricted range of even powers to convolutions of divided Bernoulli numbers and to some of their truncations \cite{MI}. These exact same $p$-residues of truncated convolutions of divided Bernoulli numbers are studied independently in $2020$ \cite{LEV3} using the $p$-adic analysis of some polynomial with $p$-adic integer coefficients whose study got initiated in \cite{LEV1}. In this study, the unsigned Stirling numbers of the first kind on $p$ letters are involved. A congruence is obtained relating these residues of truncated convolutions to the residues of some full convolutions. By confronting Miki's work and our work of \cite{LEV3}, we immediately retrieve Reijo Ernvall and Tauno Mets\"ankyla's congruence \cite{EM} (holding for the same range of powers) from $1991$, which instead of the complicated convoluted form by Miki, provides a much simpler expression for the considered residual sum in terms of the second residue in the $p$-adic Hensel expansion \cite{GOU} of a certain difference of two divided Bernoulli numbers. We show further that Ernvall and Mets\"ankyla's method does generalize to odd powers. This is the purpose of our Theorem $1$.

At this point, the $p$-adic analysis of the unsigned Stirling numbers of the first kind on $p$ letters resurfaces in order to generalize to base $a$ with $1\leq a\leq p-1$ some congruences of \cite{LEV2} that were achieved in base $2$. We provide the residue of a powers of $a$ weighted sum of divided Bernoulli numbers whose indices range between $1$ and $p-2$, for a set integer $a$ and where the powers of $a$ coincide with the indices of the divided Bernoulli numbers. This residue is as simple in expression as a sum modulo $p$ of the Wilson quotient and of the Fermat quotient in base $a$. By summing over the bases such sums after these have been multiplied by some fixed power of $a$ and by the squared Fermat quotient $q_a^2$, we thus get amongst other terms a sum of powers of integers weighted by the third power of the Fermat quotients. It is the purpose of Theorem $2$. By combining further Theorem $2$ with Theorem $1$ evoked before, we may conveniently adjust the range of indices that concerns the divided Bernoulli numbers. However, this is not yet enough to generalize the Friedmann-Tamarkine-Lehmer-Ernvall-Mets\"ankyla congruence in an appealing way, that is one which is computationally efficient. So, we continue working on the same sum and apply to it the Ernvall-Mets\"ankyla congruence. We must also generalize Ernvall and Mets\"ankyla's congruence to the unstudied cases when the power is respectively $0$ or $2$. This is the purpose of Theorems $3$ and $4$ respectively. Once this is done, we obtain Theorem $5$ which provides a satisfactory congruence modulo $p$ for the sum of powers of integers weighted by the third power of the Fermat quotient in terms of only Bernoulli numbers. Independently however, we get a much simpler form -- with a number of terms that is independent from the prime $p$ -- by using a result by Zhi-Hong Sun from \cite{SU2} on sums of powers of integers. By comparing both expressions, we obtain new congruences on Bernoulli numbers generalizing those of \cite{LEV3}. The simpler version of Theorem $5$ is the one which is retained for later use.
\\
\indent From there, we have all the ingredients to tackle the main congruence of the paper.
This time, the conjugates to the Stirling numbers, namely the multiple harmonic sums on $p-1$ integers come into play. We will shortly discuss why. In \cite{LEV4}, we had obtained by means of the Newton formulas and using some sophisticated $p$-adic expansions of the generalized harmonic numbers on $p-1$ integers modulo $p^4$ by Zhi-Hong Sun \cite{SU2} a congruence providing these sums modulo $p^4$, see Theorem $3$ of \cite{LEV4}. Sums that are related to those discussed in the previous paragraph are present in the expression. Another core term appearing in the expression is a differential term composed of two cubic convolutions of divided Bernoulli numbers. This is where Ira Gessel's identity from the title now enters the scene. The Gessel identity \cite{GE} is a generalization of the Miki identity to products of three divided Bernoulli numbers. It involves a multinomial convolution of three divided Bernoulli numbers, a binomial convolution of divided Bernoulli numbers, the harmonic numbers and the multiple harmonic sums. It plays a key role in processing the differential term. Because we may write the multiple harmonic sum $\mathcal{H}_{\lbrace s\rbrace^{2n}=1;p-1}$ in yet another way by relating it modulo $p^4$ to the Stirling numbers on $p$ letters and with $2n+1$ cycles (those are computable in an easier manner than their homologs), we are able to derive the interesting congruence modulo $p$ announced earlier. \\

Another part of the paper is devoted to studying congruences concerning sums of powers of integers weighted by the $p$-adic integer roots of some polynomials with $p$-adic integer coefficients, thus generalizing techniques arising from \cite{WJ} and \cite{EM}. We provide these congruences modulo $p^3$. The paper is built in such a way that it alternates between the different interconnected projects in order to make the computations nicer to follow. For clarity and in order to distinguish them from the present results, the past results from the other authors have been numbered with $0.i$ concerning the $i$-th Theorem or Proposition.

\subsection{Some notations}
We shall introduce some standard or some personal notations which we shall use all throughout the paper. Given a prime $p$, we denote by $w_p$ the Wilson quotient and by $q_a$ the Fermat quotient in base $a$ with $1\leq a\leq p-1$. When dealing with Bernoulli numbers, we denote by $\mb_t$ the $t-th$ divided Bernoulli number, by $\mc\mb(2k)$ a convolution of divided Bernoulli numbers of order $2k$ and by $\mbt(2k)$ a sum of products of three divided Bernoulli numbers whose indices sum up to $2k$. Also, using the same notation as in \cite{LEV3}, we denote by $\mt\mc\mb(p+1-2n, p-3)$ the following truncated convolutions of divided Bernoulli numbers:
$$\mt\mc\mb(p+1-2n,p-3):=\sum_{k=p+1-2n}^{p-3}\mb_k\mb_{2(p-1)-2n-k}$$
Throughout $\S\,2$, we also deal with binomial (resp multinomial) convolutions and thus introduce some notations for these as well. Namely,
$$b\,\mc\mb(2n):=\sum_{i=2}^{2n-2}\binom{2n}{i}\mb_i\mb_{2n-i}$$
and
$$m\,\mbt(2n):=\sum_{i+j+k=2n}\binom{2n}{i,j,k}\mb_i\mb_j\mb_k$$
By a theorem due to Clausen \cite{CL} and independently Von Staudt \cite{VS}, the denominator of a Bernoulli number $B_k$ consists of products of primes $p$ of multiplicity one, such that $p-1|k$. In particular, $pB_{p-1}=-1\mpu$. A conjecture of Agoh \cite{AG2} and independently Giuga \cite{GI} claims that $nB_{n-1}=-1\;\text{mod}\;n$ if and only if $n$ is a prime. \\
We will denote the Agoh-Giuga quotient by
$$(pB_{p-1})_1:=\frac{1+pB_{p-1}}{p}$$
When $x$ is a $p$-adic integer, we denote by $(x)_i$ its $(i+1)$-th $p$-residue in its Hensel expansion (see e.g. \cite{GOU}), that is:
$$x=\sum_{i=0}^{\infty}(x)_ip^i$$
It is important to underline that this notation does not apply to nor match with $(pB_{p-1})_1$ which got actually defined independently by being the Agoh-Giuga quotient. \\
A quite present Hensel residue throughout the paper will be:
$$\md_i:=(\mb_{p-1+i}-\mb_i)_1$$
\noindent Convolutions and truncated convolutions of these with ordinary or divided Bernoulli numbers are involved and we will thus use notations such as $\mc B\md(2n)$, $\mt\mc B\md(p+1-2n,p-3)$ and $\mc\mb\md(2n)$, $\mt\mc\mb\md(p+1-2n,p-3)$ to denote these respective entities.\\
For instance,
$$\mt\mc B\md(p+1-2n,p-3):=\sum_{i=p+1-2n}^{p-3}B_i\md_{2(p-1)-2n-i}$$
Finally, $\mh_t$ denotes the harmonic number of order $t$. Using standard notations, generalized harmonic numbers are denoted by $H_{p-1,k}$ or simply $H_k$ when it is clearly understood that the order of the sum is $p-1$. When the order of the sum is not $p-1$, we will rather denote these numbers by $\mh_{n,k}$ where the first index refers to the order of the sum and the second index refers to the index of the power. So, we have:
$$\mh_t:=\sum_{i=1}^t\frac{1}{i}\;\;\qquad H_k:=\sum_{x=1}^{p-1}\frac{1}{x^k}\;\;\qquad \mh_{n,k}:=\sum_{x=1}^n \frac{1}{x^{k}}$$
Last, regarding multiple harmonic sums, we set
$$\as_{k,n}:=\sum_{1\leq i_1<\dots< i_k\leq n}\frac{1}{i_1\dots i_k}$$
We recall that by some version of Wolstenholme's theorem, the following two congruences $H_1=0\mpd$ and $H_2=0\mpu$ hold, see \cite{WO}.

\subsection{Some historical background}
Sums of powers of integers weighted by the Fermat quotients or by powers of the Fermat quotients have drawn the interest of many mathematicians since the beginning of the twentieth century. \\

Sums of powers of integers weighted by the Fermat quotients first appear in $1909$ in \cite{FT}. Given an even integer $t$,
the authors show that:
$$\sum_{a=1}^{p-1}q_aa^t=\begin{cases}-\mb_t \mpu &\text{if $t\neq 0\;(mod\,p-1)$}\\
w_p\;\;\mpu&\text{if $t=0\;(mod\,p-1)$}
\end{cases}$$
Their proof is based on studying the number of numbers divisible by $p$ in the sequence of numbers
$$1,2,3,\dots,n$$
and noticing that this number is congruent modulo $p$ to
$$\sum_{y=1}^n(1-y^{p-1})$$
A consequence of Friedmann and Tamarkine's result is the following. \\A pair $(p,t)$ is an irregular pair if and only if
$$\sum_{a=1}^{p-1}q_aa^t=0\mpu$$
The search for irregular primes has been of importance \cite{BU}\cite{SE}\cite{WA} ever since K\"ummer \cite{KUM} showed that Fermat's last theorem holds when the exponent is a regular prime. An important feature of irregular primes is that the numerators of the divided Bernoulli numbers consist of products of powers of irregular primes, except for $\mb_n$ with $n\in\lbrace 2,4,6,8,10,14\rbrace$ when it is $1$. A study on the irregular prime divisors of the Bernoulli numbers gets pursued by Wells Johnson in \cite{JO1}.\\

Later in $1938$, Emma Lehmer generalizes the Friedmann-Tamarkine congruence \cite{LEM}.
She proves that:
$$\begin{array}{ccccc}
\sum_{a=1}^{p-1}a^{2t+1}q_a&=&-pB_{2t}&\mpd&\text{if $2t\neq 0,2\;\text{(mod $p-1$)}$}\\
\sum_{a=1}^{p-1}a^{2t}q_a&=&B_{p-1+2t}-B_{2t}&\mpd&\text{if $2t\neq 0,2\;\text{(mod $p-1$)}$}
\end{array}$$

\noindent Her proof is based on congruences concerning sums of powers of integers and on the obvious equality $a^{p-1+t}-a^t=a^tpq_a$.\\

\indent We note that the latter two congruences involve non divided Bernoulli numbers. In fact, Ernvall and Mets\"ankyla show in their $1991$ paper that the differences
$\mb_{p-1+2t}-\mb_{2t}$ where the Bernoulli numbers have now been replaced with divided Bernoulli numbers also relate modulo $p^2$ to these sums of powers of integers, but this time weighted by squared Fermat quotients. Explicitly, they generalize the K\"ummer congruences \cite{IR} by showing that:
\newtheorem*{Theorem}{Theorem $0.1$}
\begin{Theorem} Due to Ernvall and Mets\"ankyla $1991$ \cite{EM}. Let $t$ be an even integer. Then,
$$\forall 4\leq t\leq p-3,\;\;\mb_{p-1+t}=\mb_t-\frac{p}{2}\sum_{a=1}^{p-1}q_a^2a^t\;\;\mpd$$
\end{Theorem}
\noindent Their theorem follows from another of their own congruences which we recall below. First, it will be necessary to introduce the Teichm\"uller characters.
These arise for instance from factorizing the polynomial $X^{p-1}-1$ in $\mathbb{Z}_p[X]$.
Indeed, by Hensel's lemma, each integer $1\leq a\leq p-1$ lifts to a unique $p$-adic integer root
$\om(a)$ such that $|w(a)-a|_p\leq\frac{1}{p}$. The polynomial $f(X)$ factors as
$$f(X)=\prod_{a=1}^{p-1}(X-\om(a))$$
Writing for each integer $a$ with $1\leq a\leq p-1$, $$\om(a)=a+p\,v_a\qquad\text{some $a\in\mathbb{Z}_p$},$$
Ernvall and Mets\"ankyla's congruence reads as follows.
\newtheorem*{Prop}{Proposition $0.1$}
\begin{Prop} (due to Ernvall and Mets\"ankyla, \cite{EM} $1991$)
Let $t$ be an even integer with $4\leq t\leq p-3$. Then,
$$\mb_t=-\sum_{a=1}^{p-1}a^{t-1}v_a-\frac{t-1}{2}p\sum_{a=1}^{p-1}a^{t-2}v_a^2\mpd$$
\end{Prop}
\noindent Ernvall and Mets\"ankyla's proof is inspired from \cite{WJ} and is simply based on Newton's formula and sums of powers of integers. Let $t\geq 1$. Denoting by $\sigma_t$ the elementary symmetric polynomials in the roots of $f$, namely
$$\sigma_t=\sum_{1\leq a_1,\dots,a_t\leq p-1}\om(a_1)\dots\om(a_t)$$
and by $s_t$ the sums $$s_t=\sum_{a=1}^{p-1}\om(a)^t,$$ we have by Newton's formula:
$$s_t=\sigma_1s_{t-1}-\sigma_2s_{t-2}+\dots+(-1)^{t-1}t\sigma_t,\;\;\;\forall t=1\dots p-1$$
Moreover, if follows from the expansion
$$f(X)=X^{p-1}-\sigma_1X^{p-2}+\sigma_3X^{p-3}+\dots-\sigma_{p-2}X+\prod_{a=1}^{p-1}\om(a)$$
that $$\sigma_t=0,\;\;\;\forall t=1,\dots,p-2$$
Therefore, we have $s_t=0\;\forall t=1,\dots,p-2$. Substituting $w(a)=a+pv_a$ in the equation
$$\sum_{a=1}^{p-1}\om(a)^t=0,$$ expanding the $t$-th power and reducing modulo $p^3$ yields the result (see also \cite{LEV1} for the relevant developments on the sums of powers. This is where we need to exclude $t=2$, namely $\sum_{a=1}^{p-1}a^t=pB_t\;\mpd$ only when $t\neq 2$).
Because we have $\om(a)^{p-1}=1$, we also have $$\sum_{a=1}^{p-1}\om(a)^{p-1+t}=0$$
We thus get:
$$\mb_{p-1+t}-\mb_t=\sum_{a=1}^{p-1}a^{t-1}(-a^{p-1}+1)v_a+\frac{p}{2}\sum_{a=1}^{p-1}v_a^2a^{t-2}\mpd$$
Next, from expanding $(a+pv_a)^p=a+pv_a$ modulo $p^2$, we see that
$$v_a=aq_a\,\mpu$$
The theorem follows. This result is fundamentally used in \cite{LEV4} where a congruence modulo $p^4$ gets provided for $(p-1)!$ in terms of Bernoulli numbers, generalizing Wilson, Glaisher and Sun's own congruences modulo $p$, $p^2$ and $p^3$ respectively.\\

In \cite{MI}, on his way to proving his identity which has become commonly known as "Miki's identity", Hiroo Miki finds a congruence modulo $p^2$ for $$\sum_{a=1}^{p-1}a^tq_a,\;\;4\leq t\leq p-3$$ This congruence involves convolutions of divided Bernoulli numbers, binomial convolutions of divided Bernoulli numbers, truncated convolutions of divided Bernoulli numbers as well as the Agoh-Giuga quotient and harmonic numbers of order $t$.
He first finds an expansion for $aq_a,\;1\leq a\leq p-1$ to the modulus $p^2$. We state Miki's result below using our own notations.

\newtheorem*{Theo}{Theorem $0.2$}
\begin{Theo} Due to Miki $(1978)$ \cite{MI}.
Let $t$ be an even integer with $4\leq t\leq p-3$. Let $n$ be the integer such that $t=p-1-2n$. Then,
\begin{equation}\begin{split}\sum_{a=1}^{p-1}a^tq_a=&-\mb_t+\Bigg\lbrace\Big((pB_{p-1})_1+\big(1-(pB_{p-1})\big)t+\mh_t\Big)\mb_t+
\frac{t-2}{2}\mc\mb(t)\\&+\frac{1}{2}b\,\mc\mb(t)+\frac{t-1}{2}\mt\mc\mb(p+1-2n,p-3)
\Bigg\rbrace\,p\;\;(mod\,p^2)\end{split}\end{equation}
\end{Theo}

By doing a similar work on the weighted sums
$$\sum_{a=1}^{p-1}a^tq_a^2,\;\;4\leq t\leq p-3,$$
he obtains a congruence modulo $p$ for these sums in terms of convolutions of divided Bernoulli numbers, of truncated convolutions of divided Bernoulli numbers and of the Agoh-Giuga quotient. \\

Independently in \cite{LEV3}, we obtain a congruence relating residues of convolutions of divided Bernoulli numbers with residues of truncated convolutions of Bernoulli numbers. By using Ernvall and Mets\"ankyla's theorem, our Theorem $1$ point $(i)$ of \cite{LEV3} is Lemma $4$ of \cite{MI} with $m=p-1-2n$ and $6\leq m\leq p-5$.\\

Before closing this introduction, we recall below Miki's identity and its generalization by Gessel. The identity by Gessel relates a multinomial convolution of three divided Bernoulli numbers to a convolution of three divided Bernoulli numbers, while Miki's identity relates a binomial convolution of divided Bernoulli numbers to a convolution of divided Bernoulli number.
We state both results below using our own notations.

\newtheorem*{Theom}{Theorem $0.3$}
\begin{Theom} Miki's identity $1978$ \cite{MI}.\\
$$\sum_{i=2}^{m-2}\mb_i\mb_{m-i}=\sum_{i=2}^{m-2}\binom{m}{i}\mb_i\mb_{m-i}+2\mh_m\mb_m$$
\end{Theom}

\newtheorem*{Theore}{Theorem $0.4$}
\begin{Theore} Gessel's identity $2005$ \cite{GE}.\\
\begin{equation*}\begin{split}\sum_{\begin{array}{l}i+j+k=n\\i,j,k\geq 2\end{array}}\binom{n}{i,j,k}&\mb_i\mb_j\mb_k+3\mh_n\sum_{i=2}^{n-2}\binom{n}{i}\mb_i\mb_{n-i}+6\as_{2,n}\mb_n\\
&=\mbt(n)+\frac{n^2-3n+5}{4}\mb_{n-2}
\end{split}\end{equation*}
\end{Theore}

In $1982$, Shiratani and Yokoyama gave another proof for Miki's identity. Their proof of \cite{SY} is very different from Miki's original proof as it uses integration and analysis methods. It wasn't until $2005$ that Gessel gave a much simpler proof of Miki's identity and one which generalizes as in Theorem $0.4$. Gessel's proof is based on two different expressions for Stirling numbers of the second kind. An identity related to Miki's was found in $2000$ by Faber and Pandharipande \cite{FP} and proven by Zagier. Some authors like \cite{AR} have then been talking about the Miki-Zagier-Gessel identity. By using his approach, Gessel could find a generalization of the Faber-Pandharipande identity, which is also a generalization of the Miki identity. Shortly after in the same year, Crabb \cite{CR} gave a very short and simple proof for this generalization originally due to Gessel, from a functional equation for the generating function
$$\sum_{n\geq 0}\frac{B_n(\lambda)}{n!}X^n$$
Still around the same time, Schubert and Dunne found out that Miki's identity arises naturally in a certain computation in perturbative quantum field theory. They use a different generating function which plays an important role in quantum field theory computations in order to prove Miki's identity. They also prove the Faber-Pandharipande identity by using yet another different generating function. By noticing that both generating functions are related in a certain way, they also find a novel convolution identity. By their approach, they also find a cubic generalization of the Faber-Pandharipande identity. Further, they outline the generalization of the method to the derivation of convolution identities of arbitrary order. Their work appears in \cite{DS}.

\newtheorem{Remark}{Remark}
\begin{Remark}
In $1995$, in \cite{KD}, Dilcher studies multinomial convolutions of ordinary Bernoulli numbers where the number N of Bernoulli numbers in the product is greater than or equal to $5$, thus generalizing formulas by Euler (when $N=2$) and by Sitaramachandrarao and Davis \cite{SD} (when $N=3,4$).
\end{Remark}
\begin{Remark}
In $2010$ in \cite{GO}, the authors derive yet another type of quadratic convolution identities involving Bernoulli numbers, not obviously related to any included in \cite{SD}.
\end{Remark}
\begin{Remark}
In $2016$ in \cite{AG}, Agoh gave yet another proof of Miki's identity, based on Faulhaber's formula for the sums of powers of integers (sometimes known as the "Bernoulli formula").
\end{Remark}

Before stating the results of the current paper, we state below a weaker version of Miki's lemma $1$ of \cite{MI} which we shall use several times in the current paper.
\newtheorem*{Pro}{Proposition $0.2$}
\begin{Pro} Weaker version of Lemma $1$ of \cite{MI} by Miki. Let $a$ be an integer with $1\leq a\leq p-1$. Then, we have:\\
$$aq_a=-1+(1-(pB_{p-1})_1)a-\sum_{k=1}^{p-3}\frac{1}{p}\binom{p}{k}B_ka^{p-k}\mpu$$
\end{Pro}

We briefly recall a proof of this fact. Denote the sums of powers of integers by
$$S_m(n):=1^m+2^m+\dots+(n-1)^m$$
Then, Bernoulli's formula reads (see e.g. \cite{LEV1}), (in the formula below $B_1=-\frac{1}{2}$),
$$(m+1)S_m(n)=\sum_{i=1}^{m+1}\binom{m+1}{i}B_{m+1-i}n^i$$
In particular, letting $m=p-1$ and $n=a$, we have for each integer $a$ with $1\leq a\leq p-1$:
$$pS_{p-1}(a)=a^p+apB_{p-1}+\sum_{i=2}^{p-1}\binom{p}{i}B_{p-i}a^i$$
And after the change of indices $j=p-i$, we obtain:
$$pS_{p-1}(a)=a^p+apB_{p-1}+\sum_{j=1}^{p-2}\binom{p}{j}B_ja^{p-j}$$
Since $a^p=a+aq_ap$ and $pB_{p-1}=-1+p(pB_{p-1})_1$, this equality rewrites as:
$$S_{p-1}(a)=aq_a+(pB_{p-1})_1a+\sum_{j=1}^{p-2}\frac{1}{p}\binom{p}{j}B_ja^{p-j}$$
On the other hand, since $k^{p-1}=1+q_kp$, we have:
$$S_{p-1}(a)=a-1+p\sum_{k=1}^{a-1}q_k$$
Comparing the latter two identities yields Proposition $0.2$.

\subsection{New results}

Our first theorem is based on a special case of Ira Gessel's identity (listed as Theorem $0.4$ of $\S\,1.3$) and on Hiroo Miki's identity (listed as Theorem $0.3$ of $\S\,1.3$). It also relies on a result by Jianqiang Zhao from $2007$ \cite{ZO}, which combined with our work of \cite{LEV3} provides the harmonic number $\mh_{p-1}$ modulo $p^3$. Recall that this number is zero modulo $p^2$ by Wolstenholme's theorem. The proof of the theorem below is saved for the very end of the paper.

\newtheorem*{Theor}{Theorem $0$}
\begin{Theor}
\begin{equation*}
\sum_{i=2}^{p-3}\mb_i\mb_{p-1-i}\mh_{p-1-i}=2\mb_{p-3}\;\;\mpu
\end{equation*}
\end{Theor}

\newtheorem*{Coroclaire}{Corollary $0$}
\begin{Coroclaire}
\begin{equation*}\sum_{i=2}^{p-3}\mh_i\,B_i\mb_{p-1-i}=\sum_{i=2}^{p-3}\mh_i\mb_iB_{p-1-i}=-\mb_{p-3}\mpu
\end{equation*}
\end{Coroclaire}

Our next results serve as preliminary results towards the core theorem of the current paper listed below as Theorem $6$. However, these results are interesting in their own sake and generalize the works of the group of mathematicians evoked before.
By generalizing Reijo Ernvall and Tauno Mets\"ankyla's congruences described in Proposition $0.1$ and Theorem $0.1.$ of $\S\,1.3$ to odd integers $t$, we obtain a "Lehmer type" congruence like stated in Theorem $1$ below. This is a necessary step towards
finding novel congruences concerning sums of powers of integers weighted by the third power of the Fermat quotient, thus pushing the study of such sums "one Fermat quotient power" further.

\newtheorem{Th}{Theorem}

\begin{Th}
Let $t$ be an odd integer with $5\leq t\leq p-2$. Then,
\begin{equation*}
B_{p-1+t-1}-B_{t-1}=-\sum_{a=1}^{p-1}a^tq_a^2\mpu
\end{equation*}
\end{Th}

\noindent By comparing this congruence modulo $p$ to the one of Emma Lehmer modulo $p^2$, we derive in turn the following statement.

\newtheorem{Cor}{Corollary}

\begin{Cor} Let $t$ be an odd integer with $5\leq t\leq p-2$ and let $\om(a)=a+pv_a$ be the Teichm\"uller character associated to each $a$ with $1\leq a\leq p-1$. \\The following congruences hold.\\
\underline{\textit{Version $1$.}}
\begin{equation*}
\sum_{a=1}^{p-1}q_a^2a^t=-\sum_{a=1}^{p-1}q_aa^{t-1}\mpu
\end{equation*}
\underline{\textit{Version $2$.}}
\begin{equation*}
\sum_{a=1}^{p-1}q_aa^{t-1}(1+v_a)=0\;\mpu
\end{equation*}
\end{Cor}

\noindent A straightforward consequence of Theorem $1$ is also the following.

\begin{Cor} Let $t$ be an even integer with $4\leq t\leq p-3$. Then,
\begin{equation*}
\mb_t=\sum_{a=1}^{p-1}a^{t+1}q_a^2\mpu
\end{equation*}
\end{Cor}

\noindent This is the version of Theorem $1$ which we later use. The sum studied in the next theorem below is present in the congruence for the multiple harmonic sums $\as_{2n},\,2n\leq p-5$ modulo $p^4$ of Theorem $3$ of \cite{LEV4}. This sum appears to be in connection with the sum of $2n$-th powers of integers weighted by the third power of the Fermat quotient, like follows.

\begin{Th} Fix some integer $n$. Then,
\begin{equation*}
\sum_{a=1}^{p-1}\sum_{i=1}^{p-3}\frac{\mb_i}{a^i}\frac{q_a^2}{a^{2n}}=w_p\sum_{a=1}^{p-1}\frac{q_a^2}{a^{2n}}+\sum_{a=1}^{p-1}\frac{q_a^3}{a^{2n}}\mpu
\end{equation*}
\end{Th}

A joint application of Theorem $2$, Corollary $2$ applied with $t+1=p-2-2n$ implies in turn the following result.

\begin{Cor} Fix some integer $n$. Then,
\begin{equation*}
\sum_{a=1}^{p-1}\sum_{i=2}^{p-3}\frac{\mb_i}{a^i}\frac{q_a^2}{a^{2n}}=w_p\sum_{a=1}^{p-1}\frac{q_a^2}{a^{2n}}+\sum_{a=1}^{p-1}\frac{q_a^3}{a^{2n}}+\frac{1}{2}\mb_{p-3-2n}\mpu
\end{equation*}
\end{Cor}

When $2n=0$, this sum plays a central role in \cite{LEV4} in relation to finding an expansion for $(p-1)!$ to the modulus $p^4$, a generalization of Wilson, Glaisher and Sun's results to the respective moduli $p$, $p^2$ and $p^3$. \\

Some parts of the next two statements serve as a preparation to part of Theorem $5$ below.

\begin{Th}
\begin{eqnarray*}
\sum_{a=1}^{p-1}q_a^2&=&-w_p^2-\mc\mb(p-1)\qquad\qquad\qquad\qquad\qquad\;\;\;\;\;\mpu\\
&=&-\big((pB_{p-1})_1-1\big)^2-\big(2p\,\mb_{2(p-1)}-p^2\mb_{p-1}^2\big)_2\mpu\\
&&\\
&=&\frac{pB_{2(p-1)}-2pB_{p-1}+p-1}{p^2}\qquad\qquad\;\;\;\qquad\,\mpu
\end{eqnarray*}
\end{Th}

\begin{Remark} By an unpublished result of Sun \cite{SU3} which got later generalized in a published version \cite{SU1}, we have for any non-negative integer $k$,
$$pB_{k(p-1)}=-(k-1)(p-1)+kpB_{p-1}\;\mpd$$
The last congruence of Theorem $3$ offers a proof for the base case $k=2$. Sun's congruence can then be shown by induction on $k\geq 2$ using
the fact that $p^k$ divides $\sum_{a=1}^{p-1}(a^{p-1}-1)^k$. This namely implies that:
$$p\,B_{k(p-1)}=p\,\sum_{l=1}^{k-1}\binom{k}{l}B_{(p-1)(k-l)}(-1)^{l-1}+(-1)^{k-1}(p-1)\,\mpd$$
\end{Remark}

\begin{Th}
\begin{eqnarray*}
\sum_{a=1}^{p-1}q_a^2a^2&=&\frac{1}{6}w_p-\frac{1}{4}-\mt\mc\mb(4,p-3)\qquad\qquad\qquad\;\;\,\mpu\\
&=&-2(\mb_{p+1}-\mb_2)_1-\frac{1}{2}\qquad\qquad\;\;\;\qquad\qquad\mpu\\
&&\\
&=&\frac{6pB_{2p}-12pB_{p+1}+p(p-1)(2p-1)}{6p^2}\qquad\!\mpu
\end{eqnarray*}
\end{Th}

We now state the long announced generalization of the Friedmann-Tamarkin-Lehmer-Ernvall-Mets\"ankyla congruences.
\begin{Th} Let $p$ be a prime with $p\geq 7$. Set an integer $n$ such that $2\leq 2n\leq p-5$. Then, we have:
\begin{equation*}\begin{split}
\sa \frac{q_a^3}{a^{2n}}&=-2\sum_{i=2}^{p-5-2n}\mb_i(\mb_{2(p-1)-2n-i}-\mb_{p-1-2n-i})_1+2w_p(\mb_{2(p-1)-2n}-\mb_{p-1-2n})_1\\
&\qquad -\Big(1+2(\mb_{p+1}-\mb_2)_1\Big)\mb_{p-3-2n}-\Big(w_p^2+\mc\mb(p-1)\Big)\mb_{p-1-2n}\\
&\qquad -2\sum_{i=p+1-2n}^{p-3}\mb_i(\mb_{3(p-1)-2n-i}-\mb_{2(p-1)-2n-i})_1\qquad\qquad\qquad\mpu\\
&=\frac{pB_{4(p-1)-2n}-3pB_{3(p-1)-2n}+3pB_{2(p-1)-2n}-pB_{p-1-2n}}{p^3}-\mb_{p-3-2n}\\&
\qquad\qquad\qquad\qquad\qquad\qquad\qquad\qquad\qquad\qquad\qquad\qquad\qquad\qquad\mpu
\end{split}\end{equation*}
and \begin{eqnarray*}
w_p&=&(pB_{p-1})_1-1\qquad\qquad\,\mpu\\
\mc\mb(p-1)&=&(2p\mb_{2(p-1)}-p^2\mb_{p-1}^2)_2\;\mpu
\end{eqnarray*}
(Case $2n=0$)
\begin{equation*}\begin{split}\sa q_a^3&=-\big(1+2(\mb_{p+1}-\mb_2)_1\big)\mb_{p-3}+w_p\big(w_p^2+\mc\mb(p-1)\big)\\&\qquad -2\sum_{k=2}^{p-5}\mb_k(\mb_{2(p-1)-k}-\mb_{p-1-k})_1\qquad\qquad\qquad\;\;\;\,\;\,\;\;\qquad\mpu\\
&=\frac{pB_{3(p-1)}-3pB_{2(p-1)}+3pB_{p-1}-p+1}{p^3}-\mb_{p-3}\;\;\,\;\;\;\qquad\;\;\;\;\mpu
\end{split}\end{equation*}
In each case, the second congruence arises from a different perspective.
\end{Th}
\noindent In the general case we deduce a very nice congruence which is a generalization of the congruence of Theorem $1$ point $(i)$ of \cite{LEV3}.
\begin{Cor} Let $p$ be a prime with $p\geq 11$. Let $n$ be an integer such that $4\leq 2n\leq p-7$. Then,
\begin{equation*}\begin{split}
\mt\mc\mb\md(p+1-2n, p-3)&=-\mc\mb\md(p-1-2n)+w_p\big(\mb_{2(p-1)-2n}-\mb_{p-1-2n})_1\\
&-\frac{1}{2}\Big(w_p^2+\mc\mb(p-1)\Big)\mb_{p-1-2n}\\
&-\frac{pB_{4(p-1)-2n}-3pB_{3(p-1)-2n}+3pB_{2(p-1)-2n}-pB_{p-1-2n}}{2p^3}\\
&\qquad\qquad\qquad\qquad\qquad\qquad\qquad\qquad\qquad\,\;\;\;\;\mpu
\end{split}\end{equation*}
\end{Cor}

In the case when $2n=0$, we also have a nice corollary.
\begin{Cor} The following congruence holds.
\begin{equation*}\begin{split}
-\frac{1}{2}\sa\sum_{k=2}^{p-5}\frac{\mb_k}{a^k}q_a^2=&-(\mb_{p+1}-\mb_2)_1\,\mb_{p-3}+\frac{w_p}{2}\big(w_p^2+\mc\mb(p-1)\big)\\
&-\frac{pB_{3(p-1)}-3pB_{2(p-1)}+3pB_{p-1}-p+1}{2p^3}\;\;\;\mpu
\end{split}\end{equation*}
with the residues of $w_p$ and $\mc\mb(p-1)$ already provided above.
\end{Cor}
Up to a factor $p^3$, this sum was called $\mathcal{S}$ and had served as an intermediate in \cite{LEV4} for the calculation of the residue of the same sum with the divided Bernoulli numbers being replaced with the ordinary ones. The sum $\mathcal{S}$ was actually never computed as this intermediate gently vanished during the computation.

Corollary $3$ and Theorem $5$ play a key role in the proof for Theorem $6$ below. The result also relies heavily on our past work on the multiple harmonic sums $\mathcal{H}_{\lbrace s\rbrace^{2l}=1;p-1}$ modulo $p^4$ of \cite{LEV4} and on Gessel's identity of \cite{GE}. 

\begin{Th} Let $p$ be an odd prime with $p\geq 11$ and $n$ be an integer with $4\leq 2n\leq p-7$. Then,
\begin{equation*}\begin{split}\hspace{-3cm}
\mt\mc B\md=&-\mc B\md+\mc\mb\md+\frac{1}{2}\Big(\mc\mb(p-1-2n)+\mt\mc\mb(p+1-2n,p-3)\Big)_1\\
&+2n\big(\mb_{3(p-1)-2n}-2\mb_{2(p-1)-2n}+\mb_{p-1-2n}\big)_2-\big(\mb_{2(p-1)-2n}-\mb_{p-1-2n}\big)_2\\
&\hspace{-1cm}-\frac{1}{2}\Big(2(pB_{p-1})_1-1)\mb_{p-1-2n}+2\big(\mb_{2(p-1)-2n}-\mb_{p-1-2n}\big)_1\Big)_1-w_p\big(\mb_{2(p-1)-2n}-\mb_{p-1-2n}\big)_1\\
&\hspace{-2.5cm}+\Big\lbrace(2n+\frac{1}{2})\mc\mb(p-1)+w_p\big(2n-1+\frac{2n+1}{2}w_p\big)+2n\big(pB_{p-1}-p+\frac{1}{2}p^2\mb_{p-1}^2+p\mb_{p-1}-p\,\mb_{2(p-1)}\big)_2\Big\rbrace\mb_{p-1-2n}\\
&+\frac{2n+1}{2}\frac{pB_{4(p-1)-2n}-3pB_{3(p-1)-2n}+3pB_{2(p-1)-2n}-pB_{p-1-2n}}{p^3}\;\;\;\mpu
\end{split}\end{equation*}

\end{Th}

Our last theorem deals with congruences concerning sums of powers of integers weighted by the Teichm\"uller characters and by some analogs of the Teichm\"uller characters.
\begin{Th}
Let $t$ be an integer with $4\leq t\leq p-2$. \\
$(i)$ We denote by the $\om_a$'s the $(p-1)$ $p$-adic integer roots of the polynomial
$X^{p-1}-1\in\mathbb{Z}_p[X]$. The $\om_a$'s are the Teichm\"uller characters. \\ Then, we have:
$$\sum_{a=1}^{p-1}a^{t-1}\om_a=\begin{cases}p(t-1)\mb_{p-1+t}\mpt&\text{if $t$ is even.}\\p^2\big(\frac{t}{2}-1\big)B_{t-1}\mpt&\text{if $t$ is odd}\end{cases}$$
$(ii)$ We denote by the $\Om_a$'s the $(p-1)$ $p$-adic integer roots of the polynomial
$X^{p-1}+(p-1)!\in\mathbb{Z}_p[X]$. \\Then, we have:
$$\sum_{a=1}^{p-1}a^{t-1}\Om_a=\begin{cases}p(t-1)(\mb_{p-1+t}+pw_p\mb_t)\mpt&\text{if $t$ is even.}\\p^2(\frac{t}{2}-1)B_{t-1}\qquad\qquad\;\;\;\mpt&\text{if $t$ is odd}\end{cases}$$
$(iii)$ We denote by the $\gamma_a$'s the $(p-1)$ $p$-adic integer roots of the polynomial
$X^{p-1}+pB_{p-1}\in\mathbb{Z}_p[X]$. \\Then, we have:
$$\sum_{a=1}^{p-1}a^{t-1}\gamma_a=\begin{cases}p(t-1)(\mb_{p-1+t}+p(pB_{p-1})_1\mb_t)\mpt&\text{if $t$ is even.}\\p^2(\frac{t}{2}-1)B_{t-1}\qquad\qquad\qquad\;\;\;\;\;\mpt&\text{if $t$ is odd}\end{cases}$$
\end{Th}

\section{Proofs of the Theorems}

\subsection{Where we apply Gessel's identity}

This is the most technical part of the discussion. We build upon the work of \cite{LEV4} which provides expressions for the Stirling numbers $\left[\begin{array}{l}\;\;\;p\\2n+1\end{array}\right]$ (also denoted by $A_{p-1-2n}$ for the sum of products of $p-1-2n$ distinct integers chosen amongst the first $p-1$ integers) and the multiple harmonic sums $\mathcal{H}_{\lbrace s\rbrace^{2n}=1;p-1}$ (also denoted by $\as_{2n}$ for the sum of products of $2n$ distinct reciprocals of integers chosen amongst the first $p-1$ integers) modulo $p^4$. We will go straight into the technical details.
On one hand, the conjunction of Theorem $1$ and Theorem $2$ point $(i)$ of \cite{LEV4} allow to write:
\begin{equation}\begin{split}
\as_{2n}=&\frac{p^3}{6}\mbt(p-1-2n)+p\Big(1+pw_p(1+pw_p)\Big)\mb_{p-1-2n}\\&
+\frac{4(2n+1)^2+6(2n+1)+5}{24}p^3\mb_{p-3-2n}-\frac{p^2}{2}(1+pw_p)\,\mc\mb(p-1-2n)\mpq
\end{split}\end{equation}
On the other hand, we may list the main contributors from Theorem $3$ of \cite{LEV4} as follows.
\begin{equation}\begin{split}
\as_{2n}=&\frac{p^3}{6}\mbt(p-1-2n)-\frac{2n-1}{12n}p^3(\mbt(2(p-1)-2n)-\mbt(p-1-2n))\\&+
\frac{p^3}{4n}\sum_{a=1}^{p-1}\Bigg(\sum_{i=p+1-2n}^{p-5}\frac{(2n+1)\mb_i+B_i}{a^i}+\sum_{i=2}^{p-7-2n}\frac{\mb_i+B_i}{a^i}\Bigg)
\frac{q_a^2}{a^{2n}}+OT\mpq
\end{split}\end{equation}
where the other terms $OT$ must be copied from the theorem itself. \\

We will study modulo $p^4$ the differential term: $$p^3\Delta:=p^3\Big(\mbt(2(p-1)-2n)-\mbt(p-1-2n)\Big)$$
By using the Gessel identity, we may partly reduce this study to that of another simpler differential term, this time composed of multinomial cubic convolutions, namely:
$$p^3\,m\Delta:=p^3\Big(m\mbt(2(p-1)-2n)-m\mbt(p-1-2n)\Big)$$
Before we start the computation, it is worth noting that for integers $n,i,j,k$ such that $i+j+k=n$, we have:
$$\left(\begin{array}{l}n\\i,j,k\end{array}\right)
=\frac{n!}{i!j!(n-i-j)!}=\frac{n!(n-i)(n-i-1)\dots(n-i-j+1)}{i!(n-i)!j!}=\binom{n}{j}\binom{n-i}{j}$$
We will split $m\mbt(2(p-1)-2n)$ into three sums, say $S_1$, $S_2$ and $S_3$, namely:
\begin{eqnarray*}S_1&=&\binom{2(p-1)-2n}{p-1}\mb_{p-1}\sum_{j=2}^{p-3-2n}\binom{p-1-2n}{j}\mb_j\mb_{p-1-2n-j}\\
&&\\
S_2&=&\sum_{i=2}^{p-3}\binom{2(p-1)-2n}{i}\mb_i\sum_{j=2}^{2(p-1)-2n-i-2}\binom{2(p-1)-2n-i}{j}\mb_j\mb_{2(p-1)-2n-i-j}\\
&&\\
S_3&=&\sum_{i=p+1}^{2(p-1)-2n-4}\binom{2(p-1)-2n}{i}\mb_i\sum_{j=2}^{2(p-1)-2n-i-2}\binom{2(p-1)-2n-i}{j}\mb_j\mb_{2(p-1)-2n-i-j}\\
&&\\
&=&\sum_{s=4}^{p-3-2n}\binom{2(p-1)-2n}{s}\mb_{2(p-1)-2n-s}\sum_{j=2}^{s-2}\binom{s}{j}\mb_j\mb_{s-j}
\end{eqnarray*}
By Von Staudt-Clausen's theorem, the sum $S_3$ may be treated modulo $p$. When $s\leq p-3-2n$, K\"ummer's congruence applies and yields $\mb_{2(p-1)-2n-s}=\mb_{p-1-2n-s}\mpu$. Moreover, by a straightforward application of the Chu-Vandermonde identity and from other well known binomial congruences and identities (see e.g. \cite{LEV3}), we easily derive:
\begin{equation*}
\binom{2(p-1)-2n}{s}=\binom{p-1-2n}{s}+\frac{s}{2n+1}\binom{p-1-2n}{s}\mpu
\end{equation*}

And since we also have,
\begin{eqnarray*}
m\mbt(p-1-2n)&=&\sum_{i=2}^{p-5-2n}\binom{p-1-2n}{i}\mb_i\sum_{j=2}^{p-1-2n-i-2}\binom{p-1-2n-i}{j}\mb_j\mb_{p-1-2n-i-j}\\
&&\\
&=&\sum_{s=4}^{p-3-2n}\binom{p-1-2n}{s}\mb_{p-1-2n-s}\sum_{j=2}^{s-2}\binom{s}{j}\mb_j\mb_{s-j},
\end{eqnarray*}
we therefore obtain:
\begin{eqnarray*}
S_3-m\mbt(p-1-2n)&=&\frac{1}{2n+1}\sum_{s=4}^{p-3-2n}s\binom{p-1-2n}{s}\mb_{p-1-2n-s}\sum_{j=2}^{s-2}\binom{s}{j}\mb_j\mb_{s-j}\end{eqnarray*}
After a change of indices, the right hand side above rewrites as:
$$\frac{1}{2n+1}\sum_{l=2}^{p-5-2n}(p-1-2n-l)\binom{p-1-2n}{l}\mb_l\sum_{j=2}^{p-1-2n-l-2}\binom{p-1-2n-l}{j}\mb_j\mb_{p-1-2n-l-j},$$ which in turn is congruent modulo $p$ to:
$$-m\mbt(p-1-2n)-\frac{1}{2n+1}\sum_{l=2}^{p-1-2n-4}\binom{p-1-2n}{l}B_l\sum_{j=2}^{p-1-2n-l-2}\binom{p-1-2n-l}{j}\mb_j\mb_{p-1-2n-l-j}$$
The sum to the right hand side is a multinomial cubic convolution of divided Bernoulli numbers where one of the divided Bernoulli number has been  replaced with a regular Bernoulli number. By considering thrice such a sum, we see that this convolution is nothing else than: $$\frac{p-1-2n}{3}m\mbt(p-1-2n)$$
And so,
\begin{equation}
S_3-m\mbt(p-1-2n)=-\frac{2}{3}m\mbt(p-1-2n)\mpu
\end{equation}
Next, in order to tackle $S_2$ we will split it. Indeed, when $i\geq p-1-2n$, we have $p-1<i+2n+2$, then $2(p-1)-2n-i-2<p-1$. Therefore the corresponding sum may be treated modulo $p$. We thus write:
\begin{equation*}\begin{split}
S_2&=\sum_{i=2}^{p-3-2n}\binom{2(p-1)-2n}{i}\mb_i\sum_{j=2}^{2(p-1)-2n-i-2}\binom{2(p-1)-2n-i}{j}\mb_j\mb_{2(p-1)-2n-i-j}\\
&+\sum_{i=p-1-2n}^{p-3}\binom{2(p-1)-2n}{i}\mb_i\sum_{j=2}^{2(p-1)-2n-i-2}\binom{2(p-1)-2n-i}{j}\mb_j\mb_{2(p-1)-2n-i-j}
\end{split}\end{equation*}
Moreover, the sum of the second row is congruent to zero modulo $p$ as $p$ divides the binomial coefficient $\binom{2(p-1)-2n}{i}$ for this whole range of $i$.
Further, dealing now with the first row, when $i=p-3-2n$, we have $2(p-1)-2n-i-2=p-1$ and when $i<p-3-2n$, we have $2(p-1)-2n-i-2>p-1$. In the first case, by Staudt's theorem, only the extreme indices contribute to the right hand sum as $p$ divides $\binom{p+1}{j}$ for this whole range of $j$. In the second case, when $p+1-2n-i\leq j\leq p-3$, for this whole range of $j$, $p$ divides $\binom{2(p-1)-2n-i}{j}$. Then, after the adequate reductions modulo $p$ we are left with three main contributions
$$S_2=S_{2,1}+S_{2,2}+S_{2,3}$$
with: \begin{eqnarray*}S_{2,1}&=&(p+1)p\,\mb_{p-1}\mb_2\,\binom{2(p-1)-2n}{p-3-2n}\mb_{p-3-2n}\\
S_{2,2}&=&2\sum_{i=2}^{p-5-2n}\binom{2(p-1)-2n}{i}\mb_i\binom{2(p-1)-2n-i}{p-1}\mb_{p-1}\mb_{p-1-2n-i}\\
S_{2,3}&=&2\sum_{i=2}^{p-5-2n}\binom{2(p-1)-2n}{i}\mb_i\sum_{j=2}^{p-3-2n-i}\binom{2(p-1)-2n-i}{j}\mb_j\mb_{2(p-1)-2n-i-j}
\end{eqnarray*}
Regarding $S_{2,1}$, we have
\begin{eqnarray*}
\binom{2(p-1)-2n}{p-3-2n}&=&\frac{2}{p-1-2n}\binom{p-1-2n}{p-3-2n}\mpu\\
&=&-(2n+2)\qquad\qquad\qquad\;\;\mpu
\end{eqnarray*}
Therefore, \begin{equation}
S_{2,1}=-\frac{1}{6}(n+1)\mb_{p-3-2n}\;\mpu
\end{equation}
Dealing with $S_{2,2}$ relies on proving the following lemma.
\newtheorem{Lemma}{Lemma}
\begin{Lemma}Let $i$ be an integer with $0\leq i\leq p-3-2n$.
$$\binom{2(p-1)-2n-i}{p-1}\mb_{p-1}\in\mathbb{Z}_p$$ and $$\binom{2(p-1)-2n-i}{p-1}\mb_{p-1}=-\frac{1}{2n+1+i}\;\text{mod}\;p\mathbb{Z}_p$$
\end{Lemma}
\noindent \textsc{Proof of Lemma $1$.} For this range of $i$, we have $2(p-1)-2n-i>p$. Therefore,
$$\binom{2(p-1)-2n-i}{p-1}=p\;\frac{(p+p-2-2n-i)(p+p-3-2n-i)\dots (p+1)\hat{p}(p-1)!}{(p-1)!(p-1-2n-i)!}$$
Moreover, $p\mb_{p-1}=1\mpu\mathbb{Z}_p$, hence the lemma. \\

\noindent Since, for this same range of non-zero even $i$'s, we also have: $$\binom{2(p-1)-2n}{i}=\binom{p-1-2n}{i}\frac{2n+1+i}{2n+1}\;\;\;\;\;\;\;\mpu,$$
we thus get:
\begin{equation}
S_{2,2}=\frac{n+1}{6}\mb_{p-3-2n}-\frac{2}{2n+1}\,b\,\mc\mb(p-1-2n)\;\mpu
\end{equation}

Concerning $S_{2,3}$, for the ranges of even $i$'s and $j$'s that are considered, we have:
\begin{eqnarray*}
\binom{2(p-1)-2n-i}{j}&=&\binom{p-1-2n-i}{j}\frac{2n+1+i+j}{2n+1+i}\;\;\qquad\qquad\mpu
\end{eqnarray*}
It follows that:
\begin{equation*}\begin{split}
S_{2,3}=\sum_{i=2}^{p-5-2n}\binom{p-1-2n}{i}\mb_i\sum_{j=2}^{p-3-2n-i}\frac{2n+1+i+j}{2n+1}\binom{p-1-2n-i}{j}&\mb_j\mb_{p-1-2n-i-j}\\
&\mpu
\end{split}\end{equation*}
Therefore,
\begin{equation}
S_{2,3}=\frac{2}{3}\,m\mbt(p-1-2n)\;\;\;\;\;\mpu
\end{equation}

It remains to treat $S_1$. Applying Lemma $1$ with $i=0$ immediately yields:
\begin{equation}
S_1=-\frac{1}{2n+1}\;b\,\mc\mb(p-1-2n)\;\mpu
\end{equation}

By gathering congruences $(3)-(7)$, we obtain, where we also used Miki's identity and the fact that $\mh_{p-1-2n}=\mh_{2n}\mpu$ (as $\mh_{p-1}=0\mpu$ by Wolstenholme's theorem):
\begin{equation}
p^3\;m\Delta=p^3\Big(-\frac{3}{2n+1}\mc\mb(p-1-2n)+\frac{6}{2n+1}\mh_{2n}\mb_{p-1-2n}\Big)\mpq
\end{equation}
We have done only part of the work so far. It remains to study the other differential terms arising in Gessel's identity.
By Gessel's identity applied twice and after some simplifications and use of the K\"ummer congruence, we have:
\begin{equation}\begin{split}
p^3\Delta=p^3\Bigg\lbrace &m\Delta-\frac{2n+3}{2}\mb_{p-3-2n}+3\mh_{p-1-2n}\big(b\mc\mb(2(p-1)-2n)-b\mc\mb(p-1-2n)\big)\\
&+3\Big(\frac{1}{p-2n}+\frac{1}{p+1-2n}+\dots+\frac{1}{p+p-2-2n}\Big)b\mc\mb(2(p-1)-2n)\\
&+6\as_{2,2(p-1)-2n}\mb_{2(p-1)-2n}-6\as_{2,p-1-2n}\mb_{p-1-2n}\Bigg\rbrace\mpq
\end{split}
\end{equation}
We need further computations and start with the expression on the third row which we will conveniently denote by $R_3$. We decompose:
$$\as_{2,2(p-1)-2n}=\as_{2,p-1-2n}+\sum_{i=1}^{p-1-2n}\sum_{j=p-2n}^{2(p-1)-2n}\frac{1}{ij}+\sum_{p-2n\leq i<j\leq 2(p-1)-2n}\frac{1}{ij}$$
$\as_{2,2(p-1)-2n}$ is the sum of three terms, we will denote the second term by $s_2$ and the third one by $s_3$.
In $s_2$, the range of $j$ can be split into $p-2n\leq j\leq p-1$, $j=p$ and $p+1\leq j\leq p+(p-2-2n)$. When working modulo $p\mathbb{Z}_p$, we thus have:
\begin{eqnarray*}s_2&=&-\sum_{i=1}^{p-1-2n}\sum_{j=1}^{2n}\frac{1}{ij}+\frac{1}{p}\mh_{p-1-2n}+\sum_{i=1}^{p-1-2n}\sum_{j=1}^{p-2-2n}
\frac{1}{ij}\mpu
\end{eqnarray*}
Moreover, we have (note, these types of congruences get worked out later on):
\begin{equation*}
\mh_{p-1-2n}=p\mh_{2n,2}+\mh_{2n}\;\mpd
\end{equation*}
Then, we get:
$$s_2=\mh_{2n,2}+\mh_{2n}\Big(\frac{1}{p}+\frac{1}{2n+1}\Big)\;\mpu$$
We now deal with $s_3$. Regarding $s_3$, we split the ranges of $i$ and $j$ into three groups, namely
$$\begin{array}{ccc}\begin{array}{l}p-2n\\p-(2n-1)\\\;\;\vdots\\p-1\end{array}&p&\begin{array}{l}p+1\\p+2\\\;\;\vdots\\p+(p-2-2n)\end{array}\end{array}$$
We thus obtain:
\begin{equation*}s_3=\frac{1}{p}\Big(-p\,\mh_{2n,2}-\mh_{2n}-p\mh_{p-2-2n,2}+\mh_{p-2-2n}\Big)-\mh_{2n}\mh_{2n+1}+\as_{2,2n}+\as_{2,p-2-2n}\end{equation*}
\begin{equation*}\qquad\qquad\qquad\qquad\qquad\qquad\qquad\qquad\qquad\qquad\qquad\qquad\qquad\qquad\qquad\;\;\;\mpu\end{equation*}
With
\begin{eqnarray*}\mh_{p-2-2n}&=&p\,\mh_{2n+1,2}+\mh_{2n+1}\;\mpd\\
\mh_{p-2-2n,2}&=&-\mh_{2n+1,2}\;\;\;\;\qquad\;\;\;\;\mpu\;,\end{eqnarray*}
$s_3$ then simplifies to
$$s_3=\frac{1}{p}\Big(\frac{1}{2n+1}+\frac{p}{(2n+1)^2}+p\mh_{2n+1,2}\Big)-\mh_{2n}\mh_{2n+1}+\as_{2,2n}+\as_{2,p-2-2n}\mpu$$
Moreover, we have:
\begin{eqnarray*}\as_{2,p-2-2n}&=&\as_{2,p-1-2n}+\frac{\mh_{2n+1}}{2n+1}\;\mpu\end{eqnarray*}
We now need the following intermediate result which gets listed as a lemma.
\begin{Lemma}
$$\as_{2,p-1-2n}=-\as_{2,2n}+\mh_{2n}^2\;\;\mpu$$
\end{Lemma}
\textsc{Proof of Lemma $2$.} We have:
\begin{eqnarray*}
\as_{2,p-1-2n}&=&\as_2-\sum_{p-2n\leq i<j\leq p-1}\frac{1}{ij}-\sum_{i=1}^{p-1-2n}\sum_{j=p-2n}^{p-1}\frac{1}{ij}\\
&=&-\as_{2,2n}+\sum_{i=1}^{p-1-2n}\sum_{j=1}^{2n}\frac{1}{ij}\mpu\\
&=&-\as_{2,2n}+\mh_{2n}\mh_{p-1-2n}\mpu\\
&=&-\as_{2,2n}+\mh_{2n}^2\;\;\qquad\;\;\;\;\,\mpu
\end{eqnarray*}
\hfill$\square$\\
Then, $$\as_{2,p-2-2n}=-\as_{2,2n}+\mh_{2n}\mh_{2n+1}+\frac{1}{(2n+1)^2}\;\;\mpu$$
Thus, \begin{eqnarray*}s_3&=&\frac{1}{2n+1}\Big(\frac{1}{p}+\frac{2}{2n+1}\Big)+\mh_{2n+1,2}\;\;\;\;\qquad\mpu\end{eqnarray*}
In the end, we have:
\newtheorem{Proposition}{Proposition}
\begin{Proposition}
$$\as_{2,2(p-1)-2n}=\as_{2,p-1-2n}+2\mh_{2n+1,2}+\mh_{2n+1}\Big(\frac{1}{p}+\frac{1}{2n+1}\Big)\;\mpu$$
\end{Proposition}
It follows that
\begin{equation}
p^3R_3=6\Big(\frac{\mh_{2n+1}}{2n+1}+2\mh_{2n+1,2}\Big)p^3\mb_{p-1-2n}+6\mh_{2n+1}p^2\mb_{2(p-1)-2n}\;\mpq
\end{equation}

We further deal with $R_2$. Let $u$ denote the $p$-adic integer:
$$u:=\frac{1}{p-2n}+\dots+\frac{1}{p-1}+\frac{1}{p+1}+\dots+\frac{1}{p+p-2-2n}$$
We have:
$$R_2=\frac{3}{p}\,b\,\mc\mb(2(p-1)-2n)+3\,u\,b\,\mc\mb(2(p-1)-2n)$$
We decompose:
\begin{equation*}\begin{split}
b\,\mc\mb(2(p-1)-2n)=&2\binom{2(p-1)-2n}{p-1}\mb_{p-1}\mb_{2(p-1)-2n}\\&+2\sum_{i=2}^{p-3-2n}\binom{2(p-1)-2n}{i}\mb_i\mb_{2(p-1)-2n-i}\\
&+\sum_{i=p+1-2n}^{p-3}\binom{2(p-1)-2n}{i}\mb_i\mb_{2(p-1)-2n-i}\end{split}\end{equation*}
Then,
\begin{equation*}\begin{split}p^3R_2=&3p^2\,b\,\mc\mb(2(p-1)-2n)+6\binom{2(p-1)-2n}{p-1}p\,\mb_{p-1}\mb_{p-1-2n}p^2\,u\\
&+\frac{6}{2n+1}p^3u\sum_{i=2}^{p-3-2n}\binom{p-1-2n}{i}\mb_i\mb_{p-1-2n-i}(2n+1+i)\\&+3p^3u\sum_{i=p+1-2n}^{p-3}\binom{2(p-1)-2n}{i}\mb_i\mb_{2(p-1)-2n-i}
\mpq\end{split}\end{equation*}
When $p+1-2n\leq i\leq p-3$, the binomial coefficient in the last sum is divisible by $p$, hence the latter sum simply vanishes. Further, after scrutiny, $b\,\mc\mb(2(p-1)-2n)$ is wanted modulo $p^2$ while $u$ is simply desired modulo $p$.
We have:
\begin{eqnarray*}
u&=&\mh_{p-2-2n}-\mh_{2n}\;\;\;\qquad\;\;\;\;\mpu\\
&=&H_1-\sum_{k=p-1-2n}^{p-1}\frac{1}{k}-\mh_{2n}\;\mpu\\
&=&\frac{1}{2n+1}\;\;\;\;\;\;\qquad\qquad\qquad\mpu\end{eqnarray*}
We now rewrite:
\begin{equation*}\begin{split}p^3R_2=3p^2\,b\,\mc\mb(2(p-1)-2n)-\frac{6}{(2n+1)^2}p^3\,\mb_{p-1-2n}+\frac{3}{2n+1}p^3\,b\,&\mc\mb(p-1-2n)\\
&\qquad\;\;\;\mpq\end{split}\end{equation*}
We move on to studying the binomial convolution $b\,\mc\mb(2(p-1)-2n)\,mod\,p^2\mathbb{Z}_p$.
We have:
\begin{eqnarray*}
3p^2\,b\,\mc\mb(2(p-1)-2n)&=&3p^2\mc\mb(2(p-1)-2n)-6p^2\mh_{2(p-1)-2n}\mb_{2(p-1)-2n}\\
&=&3p^2\mc\mb(2(p-1)-2n)-6p^2\mh_{p-1-2n}\mb_{2(p-1)-2n}\\&&-6p\,\mb_{2(p-1)-2n}-6p^2\,u\,\mb_{2(p-1)-2n}
\end{eqnarray*}

The first equality above holds by Miki's identity. This time, we need to know $u$ modulo $p^2$ instead of simply modulo $p$. This is routine calculation which we expand right below.
We have:
\begin{equation*}\begin{split}
\frac{1}{p-2n}&+\frac{1}{p+1-2n}+\dots+\frac{1}{p-1}\\&=-\frac{1}{(2n)^2}(p+2n)-\frac{1}{(2n-1)^2}(p+2n-1)-\dots-\frac{1}{1^2}(p+1)\qquad\mpd\\
&=-p\,\mh_{2n,2}-\mh_{2n}\qquad\qquad\qquad\qquad\qquad\qquad\qquad\qquad\qquad\qquad\;\;\,\mpd
\end{split}
\end{equation*}
\begin{equation*}\begin{split}
\frac{1}{p+1}&+\frac{1}{p+2}+\dots+\frac{1}{p+p-2n-2}\\&=-\frac{1}{1^2}(p-1)-\frac{1}{2^2}(p-2)-\dots-\frac{1}{(p-2n-2)^2}(p-(p-2n-2))\;\mpd\\
&=-p\,\mh_{p-2n-2,2}+\mh_{p-2n-2}\qquad\qquad\qquad\qquad\qquad\qquad\qquad\qquad\;\;\;\,\mpd\\
&=-pH_2+p\Big(\frac{1}{(p-1-2n)^2}+\dots+\frac{1}{(p-1)^2}\Big)\\&\;\;\;+H_1-\Big(\frac{1}{p-1-2n}+\frac{1}{p-2n}+\dots+\frac{1}{p-1}\Big)\qquad
\qquad\qquad\qquad\;\;\mpd\\
&=2p\;\mh_{2n+1,2}+\mh_{2n+1}\qquad\qquad\qquad\qquad\qquad\qquad\qquad\qquad\qquad\;\;\;\,\,\,\mpd
\end{split}\end{equation*}
Adding both contributions yields:
$$u=\frac{1}{2n+1}+p\Big(\frac{1}{(2n+1)^2}+\mh_{2n+1,2}\Big)\mod\,p^2\mathbb{Z}_p$$
We also have:
$$\mh_{p-1-2n}=p\,\mh_{2n,2}+\mh_{2n}\,\qquad\qquad\qquad\mpd\mathbb{Z}_p$$
Therefore,
\begin{equation*}\begin{split}
3p^2\,b\,\mc\mb(2(p-1)-2n)&=3p^2\mc\mb(2(p-1)-2n)\\&\;\;\;\;-12p^3\mh_{2n+1,2}\mb_{p-1-2n}-6p(1+p\mh_{2n+1})\mb_{2(p-1)-2n}\\
&=6p^2\sum_{i=2}^{p-3-2n}\mb_i\mb_{2(p-1)-2n-i}+6p\mb_{p-1}\,p\mb_{p-1-2n}\\&\;\;\;\; +3p^2\mt\mc\mb(p+1-2n,p-3)\\&\;\;\;\;
-12p^3\mh_{2n+1,2}\mb_{p-1-2n}-6p(1+p\mh_{2n+1})\mb_{2(p-1)-2n}\\&\;\;\;\qquad\qquad\qquad\qquad\qquad\qquad\qquad\qquad\qquad\qquad\mpq\\
\end{split}\end{equation*}
In order to apply Ernvall and Mets\"ankyla's congruence in the first sum above, we must impose $i\leq p-5-2n$. Hence we split this sum. 
After grouping the different terms, we obtain:

\begin{equation*}\begin{split}
p^3R_2=&\frac{3}{2n+1}p^3\mc\mb(p-1-2n)+3p^2\mc\mb(p-1-2n)\\&+3p^2\Big(\mc\mb(p-1-2n)+\mt\mc\mb(p+1-2n,p-3)\Big)-3p^3\sum_{i=2}^{p-5-2n}\sum_{a=1}^{p-1}\frac{q_a^2}{a^{2n}}\frac{\mb_i}{a^i}\\
&+6p\Big[p\mb_{p-1}-p^2\Big(2\mh_{2n+1,2}+\frac{\mh_{2n+1}}{2n+1}\Big)\Big]\mb_{p-1-2n}+6p^3(\mb_{p+1}-\mb_2)_1\mb_{p-3-2n}\\
&-6p(1+p\mh_{2n+1})\mb_{2(p-1)-2n}\qquad\qquad\qquad\qquad\qquad\qquad\qquad\qquad\mpq
\end{split}\end{equation*}
Moreover, by Theorem $1$ point $(i)$ of \cite{LEV3},
$$\mt\mc\mb(p+1-2n,p-3)+\mc\mb(p-1-2n)=-\sum_{a=1}^{p-1}\frac{q_a^2}{a^{2n}}+2\Big((pB_{p-1})_1-1\Big)\mb_{p-1-2n}\mpu$$
Then,
\begin{equation}\begin{split}
p^3R_2=&\frac{3}{2n+1}p^3\mc\mb(p-1-2n)+3p^2\mc\mb(p-1-2n)\\&+3p^3\Big(\mc\mb(p-1-2n)+\mt\mc\mb(p+1-2n,p-3)\Big)_1\\
&+3p^2\Big(2\big((pB_{p-1})_1-1\big)\mb_{p-1-2n}-\sum_{a=1}^{p-1}\frac{q_a^2}{a^{2n}}\Big)_0-3p^3\sum_{i=2}^{p-5-2n}\sum_{a=1}^{p-1}\frac{q_a^2}{a^{2n}}\frac{\mb_i}{a^i}\\
&+6p\Big[p\mb_{p-1}-p^2\Big(2\mh_{2n+1,2}+\frac{\mh_{2n+1}}{2n+1}\Big)\Big]\mb_{p-1-2n}+6p^3(\mb_{p+1}-\mb_2)_1\mb_{p-3-2n}\\
&-6p(1+p\mh_{2n+1})\mb_{2(p-1)-2n}\qquad\qquad\qquad\qquad\qquad\qquad\qquad\qquad\mpq
\end{split}\end{equation}
It remains to tackle the third term of the first row of $(10)$. We have:
\begin{equation*}\begin{split}p^3R_1=&6\mh_{p-1-2n}p\mb_{p-1}p^2\mb_{p-1-2n}\binom{2(p-1)-2n}{p-1}\\
&+3\mh_{p-1-2n}p^3\sum_{i=2}^{p-3-2n}\Big(2\binom{2(p-1)-2n}{i}-\binom{p-1-2n}{i}\Big)\mb_i\mb_{p-1-2n-i}\\&
\qquad\qquad\qquad\qquad\qquad\qquad\qquad\qquad\qquad\qquad\qquad\qquad\qquad\qquad\mpq
\end{split}\end{equation*}
As seen several times before, for the range of $i$ that is considered,
$$\binom{2(p-1)-2n}{i}=\binom{p-1-2n}{i}\frac{2n+1+i}{2n+1}\mpu$$
Then,
$$2\binom{2(p-1)-2n}{i}=\binom{p-1-2n}{i}\Big(2+\frac{2i}{2n+1}\Big)\mpu$$
It immediately follows that the difference above vanishes. 
Further, by Lemma $1$, we know that:
$$\binom{2(p-1)-2n}{p-1}\mb_{p-1}=-\frac{1}{2n+1}\mpu$$
Then, we have:
\begin{equation}
p^3R_1=-6p^3\frac{\mh_{2n}}{2n+1}\mb_{p-1-2n}\mpq
\end{equation}
By using Result $12$ of \cite{LEV3} due to Sun \cite{SU2} with $k=2$ and $b=p-1-2n$, as well as Result $10$ point $(i)$ of \cite{LEV3}, the congruence modulo $p^3\mathbb{Z}_p$ arising from $(2)$ and $(3)$ simply reads:
$$p^3\Delta=3p^2\mc\mb(p-1-2n)\mpt\mathbb{Z}_p$$
It is then routine verification that the congruence modulo $p^3\mathbb{Z}_p$ is just trivial. We thus focus our attention on the modulus $p^4\mathbb{Z}_p$. But before moving any further, this is an adequate time to prove the first few theorems of the introduction, as well as their related corollaries, culminating in the proof of Theorem $6$.

\subsection{Proofs of Theorems $1-5$}
Because the sum of importance in the previous $\S\,2.1$ is the one contained in the left hand side of the congruence of Theorem $2$, we start with the proof of Theorem $2$. \\The proof is inspired from \cite{LEV2} and uses the Stirling numbers. The unsigned Stirling number of the first kind $\left[\begin{array}{l}p\\s\end{array}\right]$ is the unsigned coefficient of $x^s$ in the falling factorial
$$x(x-1)(x-2)\dots (x-(p-1))=\sum_{s=1}^p(-1)^{s-1}\left[\begin{array}{l}p\\s\end{array}\right]x^s$$
By specializing $x=-a$ in the equality above, we obtain:
$$-a(-a-1)(-a-2)\dots(-a-(p-1))=\sum_{s=1}^{p}(-1)^{s-1}\left[\begin{array}{l}p\\s\end{array}\right](-a)^s$$
Then,
$$\frac{(p-1+a)!}{(a-1)!}=\sum_{s=1}^p\left[\begin{array}{l}p\\s\end{array}\right]a^s$$
Moreover,
\begin{eqnarray*}
(p-1+a)!&=&(p-1)!p\dots(p-1+a)\\
&=&-p(a-1)!\mpd
\end{eqnarray*}
It follows that:
$$-p=\sum_{s=1}^{p}\left[\begin{array}{l}p\\s\end{array}\right]a^s\mpd$$
By using Corollary $2$ of \cite{LEV1} originally due to Glaisher \cite{GL}, we get:
$$-p=(pB_{p-1}-p)a+\sum_{s=3}^{p-2}\frac{p}{s}B_{p-s}a^s-\frac{p}{2}a^{p-1}+a^p\;\mpd$$
It follows that:
\begin{Proposition}
$$\sum_{i=1}^{p-3}\frac{\mb_i}{a^i}=w_p+q_a\;\mpu$$
\end{Proposition}
From there, Theorem $2$ follows. A rewriting of the theorem is the following.
$$-\frac{1}{2}\sum_{a=1}^{p-1}\frac{q_a^2}{a^{2n+1}}+\sum_{a=1}^{p-1}\sum_{i=2}^{p-3}\frac{\mb_i}{a^i}\frac{q_a^2}{a^{2n}}=w_p\sum_{a=1}^{p-1}\frac{q_a^2}{a^{2n}}+\sum_{a=1}^{p-1}\frac{q_a^3}{a^{2n}}\mpu$$
We now adapt Ernvall and Mets\"ankyla's proof for the results of $\S\,1.3$ when $t$ is odd. First and foremost, we recall from \cite{SU2} that
\begin{equation}S_t=\sum_{a=1}^{p-1}a^t=pB_t+\frac{p^2}{2}tB_{t-1}+\frac{p^3}{6}t(t-1)B_{t-2}\;\mpt\end{equation}
Suppose $t$ is odd and $t\geq 3$. Then $B_t=0$. If we impose $t\neq 3$, then
$$S_t=\frac{p^2}{2}tB_{t-1}\mpt$$
Again, we have, using the same notations as in $\S\,1$:
$$(a+pv_a)^t=a^t+tpv_aa^{t-1}+p^2v_a^2a^{t-2}\frac{t(t-1)}{2}\mpt$$
Like before, when $t$ is odd and $3\leq t\leq p-2$, we have
$$\sum_{a=1}^{p-1}(a+pv_a)^t=0$$
When excluding $t=3$, reducing modulo $p^3$ yields:
\begin{equation}\frac{p^2}{2}tB_{t-1}=-tp\sum_{a=1}^{p-1}v_aa^{t-1}-p^2\frac{t(t-1)}{2}\sum_{a=1}^{p-1}v_a^2a^{t-2}\mpt\end{equation}
We thus have:
$$\frac{p}{2}B_{t-1}=-\sum_{a=1}^{p-1}v_aa^{t-1}-p\frac{t-1}{2}\sum_{a=1}^{p-1}v_a^2a^{t-2}\mpd$$
And since we also have $$\sum_{a=1}^{p-1}(a+pv_a)^{p-1+t}=0,$$
similar developments and reductions lead to:
$$\frac{p}{2}B_{p-1+t-1}=-\sum_{a=1}^{p-1}v_aa^{p-1+t-1}-p\frac{t-2}{2}\sum_{a=1}^{p-1}v_a^2a^{t-2}\mpd$$
Then, using also $v_a=aq_a\mpu$, we obtain Theorem $1$.\\
It follows from Theorem $2$ and Corollary $2$ with $t+1=p-2-2n$ that
\begin{equation}
\sum_{a=1}^{p-1}\sum_{i=2}^{p-3}\frac{\mb_i}{a^i}\frac{q_a^2}{a^{2n}}=w_p\sum_{a=1}^{p-1}\frac{q_a^2}{a^{2n}}+\sum_{a=1}^{p-1}\frac{q_a^3}{a^{2n}}+\frac{1}{2}\mb_{p-3-2n}\mpu
\end{equation}

\noindent This result can also be derived from a weaker version of Lemma $1$ of \cite{MI}. Namely, as part of Miki's congruence stated in his Lemma $1$ of \cite{MI}, we have:
$$aq_a=-1+(1-(pB_{p-1})_1)a-\sum_{k=1}^{p-3}\frac{1}{p}\binom{p}{k}B_ka^{p-k}\mpu$$
It comes:
$$\sum_{a=1}^{p-1}aq_a\frac{q_a^2}{a^{2n+1}}=-\sum_{a=1}^{p-1}\frac{q_a^2}{a^{2n+1}}+(1-(pB_{p-1})_1)\sum_{a=1}^{p-1}\frac{q_a^2}{a^{2n}}-\sum_{k=1}^{p-3}\frac{1}{p}\binom{p}{k}B_k\sum_{a=1}^{p-1}q_a^2a^{p-1-2n-k}\mpu$$
From there, by using again Corollary $2$ with $t+1=p-2-2n$, we retrieve $(16)$.\\

After several adequate applications of the Ernvall-Mets\"ankyla congruence in $(16)$, we further get:
\begin{equation}\begin{split}
\sa \frac{q_a^3}{a^{2n}}=&-2\sum_{i=2}^{p-5-2n}\mb_i(\mb_{2(p-1)-2n-i}-\mb_{p-1-2n-i})_1\\
&+\Bigg(\sa q_a^2a^2-\frac{1}{2}\Bigg)\mb_{p-3-2n}+\Bigg(\sa q_a^2\Bigg)\mb_{p-1-2n}+2w_p(\mb_{2(p-1)-2n}-\mb_{p-1-2n})_1\\
&-2\sum_{i=p+1-2n}^{p-3}\mb_i(\mb_{3(p-1)-2n-i}-\mb_{2(p-1)-2n-i})_1\mpu
\end{split}\end{equation}

In order to conclude to Theorem $5$, we first need to prove Theorems $3$ and $4$ of $\S\,1.3$. These theorems are built out of two perspectives. The first perspective relies on Miki's result providing the residue of $aq_a$ for $a$ integer with $1\leq a\leq p-1$. The second perspective relies on Sun's congruence for the sums of powers of integers issued from \cite{SU2}. We first discuss the first perspective. It treats the first two congruences in each theorem. \\

By a reformulation of Proposition $0.2$ of $\S\,1.3$ due to Miki, we have:
$$aq_a=-\frac{1}{2}+(1-(pB_{p-1})_1)a+\sum_{k=2}^{p-3}\mb_ka^{p-k}\mpu$$
Then,
$$\sum_{a=1}^{p-1}q_a^2=\sum_{a=1}^{p-1}aq_a\frac{q_a}{a}=-\frac{1}{2}\sum_{a=1}^{p-1}q_aa^{p-2}+(1-(pB_{p-1})_1w_p+\sum_{k=2}^{p-3}\sum_{a=1}^{p-1}\frac{\mb_k}{a^k}q_a\mpu$$
Moreover, Lehmer's congruence for odd powers implies that
$$\sum_{a=1}^{p-1}q_aa^{p-2}=-pB_{p-3}\mpd$$
Then, the residue of the first sum is simply zero.\\ Further, by the original Friedmann-Tamarkine congruence, we have:
$$\sa q_aa^{p-1-k}=-\mb_{p-1-k}\;\mpu$$
Then, up to sign, the last sum is the convolution $\mc\mb(p-1)$. Hence the first congruence of Theorem $3$. The second congruence follows from \cite{LEV3} which provides the residue of $\mc\mb(p-1)$ (see Result $11$ of \cite{LEV4}).\\

As for Theorem $4$, we have:
\begin{equation*}\begin{split}\sa q_a^2a^2=\sa (aq_a)(aq_a)=-\frac{1}{2}\sa q_aa+(1-(pB_{p-1})_1)\sa q_aa^2+\sum_{k=2}^{p-3}\sa& \frac{\mb_k}{a^k}q_aa^2\\
&\mpu\end{split}\end{equation*}

\noindent When $4\leq k\leq p-3$, the Friedmann-Tamarkine congruence applies and yields:
$$\sa q_aa^{p+1-k}=-\mb_{p+1-k}\;\mpu$$
Moreover, the Friedmann-Tamarkine congruence also applies to the first two terms. We thus get the first congruence of Theorem $4$. The second congruence will follow immediately from forthcoming Lemma $3$ of $\S\,2.3$\\

The proof for the last congruence in each theorem relies on Congruence $(14)$ by Z-H. Sun in \cite{SU2} and is left to the reader.

From there, Theorem $5$ follows.

\subsection{Proof of Theorem $6$}
This part closes the proof for Theorem $6$. The rest of the computations will be based on the following series of lemmas. 

\begin{Lemma}
\begin{equation*}
\mt\mc\mb(4,p-3)=2(\mb_{p+1}-\mb_2)_1+\frac{1}{12}+\frac{1}{6}(pB_{p-1})_1\;\mpu
\end{equation*}
\end{Lemma}
\textsc{Proof.} Immediately follows from Theorem $1$ point $(ii)$ of \cite{LEV3}.

\begin{Lemma}
\begin{equation*}
\mt\mc\mb(6,p-3)=\frac{7}{720}+2(\mb_{p+3}-\mb_4)_1-\frac{1}{60}(pB_{p-1})_1\;\mpu
\end{equation*}
\end{Lemma}
\textsc{Proof.} Again, this is an immediate consequence of Theorem $1$ point $(iii)$ of \cite{LEV3}.

\begin{Lemma}
\begin{eqnarray*}
\sa\sum_{i=2}^{p-5-2n}\frac{B_i}{a^i}\frac{q_a^2}{a^{2n}}&=&-2\sum_{i=2}^{p-5-2n}B_i(\mb_{2(p-1)-2n-i}-\mb_{p-1-2n-i})_1\qquad\mpu\\
\sa\sum_{i=p+1-2n}^{p-3}\frac{B_i}{a^i}\frac{q_a^2}{a^{2n}}&=&-2\sum_{i=p+1-2n}^{p-3}B_i(\mb_{3(p-1)-2n-i}-\mb_{2(p-1)-2n-i})_1\mpu
\end{eqnarray*}
\end{Lemma}
\textsc{Proof.} For the first range of $i$, we have $4\leq p-1-2n-i\leq p-3$. For the second range of $i$, we have $4\leq 2(p-1)-2n-i\leq p-3$.
In both cases the Ernvall-Mets\"ankyla congruence applies and yields the result.
\begin{Lemma}
\begin{equation*}
\sa\sum_{i=2}^{p-5-2n}\frac{\mb_i}{a^i}\frac{q_a^2}{a^{2n}}=-2\sum_{i=2}^{p-5-2n}\mb_i(\mb_{2(p-1)-2n-i}-\mb_{p-1-2n-i})_1\mpu
\end{equation*}
\end{Lemma}
\textsc{Proof.} This range for $i$ satisfies to the conditions of application of the Ernvall Mets\"ankyla's congruence.

More computational efforts using Lemma $3-6$ lead to Theorem $6$ which is thus entirely proven. The next part deals again with the techniques exposed in the introduction.

\subsection{Proof of Theorem $7$}
In the case when $t$ is even, point $(i)$ of Theorem $7$ follows from a joint application of Proposition $0.1$ and Theorem $0.1$ of $\S\,1.3$, both due to Ernvall and Mets\"ankyla.\\
In the case when $t$ is odd, the sum of powers $S_t$ verifies
\begin{equation}S_t=\frac{p^2}{2}tB_{t-1}\mpt\end{equation}
And we have seen in $\S\,2.2$ (cf. Congruence $(15)$) that
$$p\sum_{a=1}^{p-1}a^{t-1}v_a=-\frac{p^2}{2}B_{t-1}-\frac{p^2}{2}(t-1)\sum_{a=1}^{p-1}a^tq_a^2\;\;\mpt$$
By summing both contributions and using our Corollary $2$, we obtain the result of $(i)$ which is thus entirely proven. \\
Regarding point $(ii)$, we will study the polynomial $$X^{p-1}+(p-1)!\in\mathbb{Z}_p[X]$$
in the same way as Ernvall and Mets\"ankyla had used the polynomial $$X^{p-1}-1\in\mathbb{Z}_p[X]$$ in order to prove their congruences. The first polynomial had played a key role in \cite{LEV1} in order to recover Glaisher's result that $(p-1)!=pB_{p-1}-p\mpd$ and in order to find a formula for $(p-1)!$ modulo $p^3$ by a different method than the one of Sun \cite{SU2}, see Corollary $6$ in $\S\,1$ of \cite{LEV1}. The analogs of the Teichm\"uller characters are the $(p-1)$ $p$-adic integer roots of $X^{p-1}+(p-1)!$ which we denote by $\Om_a$, $1\leq a\leq p-1$. Moreover, we define the $w_a$'s such that:
$$\Om_a=a+pw_a$$
It is shown in \cite{LEV1} that
$$pw_a=a(1+(p-1)!+pq_a)\mpd$$
This is the purpose of Lemma $1$ in $\S\,2$ of \cite{LEV1}. Thus, we have:
$$w_a=a(w_p+q_a)\;\mpu$$
By the same argument as in $\S\,1$, we have:
$$\sum_{a=1}^{p-1}(a+pw_a)^t=0$$
We assume $t\neq 2$. Expanding modulo $p^3$ leads to an analog of the Ernvall and Mets\"ankyla's congruence, namely
$$\mb_t=-\sum_{a=1}^{p-1}w_aa^{t-1}-p\frac{t-1}{2}\sum_{a=1}^{p-1}w_a^2a^{t-2}\mpd,$$
where $v_a$ was replaced with $w_a$.
By substituting:
\begin{eqnarray*}
w_a^2&=&a^2w_p^2+q_a^2a^2+2a^2q_aw_p\mpu,
\end{eqnarray*}
we then obtain for even $t$:
$$\mb_t+\sum_{a=1}^{p-1}w_aa^{t-1}=-p\frac{t-1}{2}\sum_{a=1}^{p-1}q_a^2a^t-p(t-1)w_p\sum_{a=1}^{p-1}q_aa^t\;\mpd$$
By the Lehmer congruence applied with even powers modulo $p$ and by the Ernvall-Mets\"ankyla's congruence, we further obtain:
\begin{equation}
\sum_{a=1}^{p-1}w_aa^{t-1}=(t-1)\mb_{p-1+t}-t\mb_t+p(t-1)w_p\mb_t\;\mpd
\end{equation}
The case $t$ is even of Theorem $7$ point $(ii)$ follows. When $t$ is odd, the congruence reads instead:
$$\sum_{a=1}^{p-1}w_aa^{t-1}=-\frac{p}{2}B_{t-1}-p\frac{t-1}{2}\sum_{a=1}^{p-1}q_a^2a^t-p(t-1)w_p\sum_{a=1}^{p-1}q_aa^t\mpd$$
The first sum to the right hand side of the congruence is known modulo $p$ by our Corollary $2$. As for the second sum it is known modulo $p^2$ by Lehmer's result applied with odd powers. The latter sum thus vanishes from the congruence since $t\neq 1\,\text{mod}\,(p-1)$. We get:
\begin{equation}\sum_{a=1}^{p-1}w_aa^{t-1}=-pB_{t-1}\;\mpd\end{equation}
We deduce point $(ii)$ of Theorem $7$ in the case when $t$ is odd. \\

We now present a different proof for Congruences $(19)$ and $(20)$, one which does not refer to Ernvall and Mets\"ankyla's original proof of \cite{EM}.
We use an expansion of $w_a$ one $p$-power further. This is achieved in \cite{LEV1}. Our $(w_a)_1$ is the $t_a^{(1)}$ of \cite{LEV1}. Also, $\delta_i(a)$ of \cite{LEV1} is $(q_a)_i$. Using the current notations, Lemma $3$ of \cite{LEV1} asserts that:
$$(w_a)_1=a\Big((q_a)_0+(q_a)_1+\Bigg(\sum_{i=1}^{p-1}(q_i)_0\Bigg)^2+(1+(q_a)_0)\,\sum_{i=1}^{p-1}(q_i)_0\Big)\;\;\text{mod}\,p$$
Recall also from Theorem $1$ of \cite{LEV1} that
$$\sum_{a=1}^{p-1}q_a=w_p\;\mpu$$
Then, using again $$\mb_t=-\sum_{a=1}^{p-1}a^tq_a\mpu,$$ we have:
$$\sum_{a=1}^{p-1}w_aa^{t-1}=pw_p(t-1)\mb_t-p\mb_t+\sum_{a=1}^{p-1}a^tq_a\mpd$$
But by the Lehmer congruence applied with $t$ even,
$$\sum_{a=1}^{p-1}a^tq_a=p\mb_t+(t-1)\mb_{p-1+t}-t\mb_t\mpd$$
Combining both congruences yields $(19)$. \\
When $t$ is odd, we rather apply the Lehmer congruence with odd powers and it leads to $(20)$.\\

It remains to deal with $(iii)$. We define in turn the $p$-adic integers $z_a$'s such that for each integer $a$ with $1\leq a\leq p-1$,
$$\gamma_a=a+pz_a$$
By Hensel's lifting algorithm, the residue $z_a$ gets determined by
$$(a+pz_a)^{p-1}+pB_{p-1}\in p^2\mathbb{Z}_p$$
We thus have:
$$z_a=a\big(q_a+(pB_{p-1})_1\big)\mpu$$
Then the sketch of the proof with respect to the Ernvall-Mets\"ankyla type of proof is identical as in $(ii)$ with $w_p$ being replaced with $(pB_{p-1})_1$. It yields the result of $(iii)$ straight away. Again, another proof consists of pushing the expansion for $z_a$ one $p$-power further as in Proposition $3$ below whose proof relies on Hensel's lifting algorithm.
\begin{Proposition}
$$z_a=a(q_a+(pB_{p-1})_1)+ap\big(1+(pB_{p-1})_1\big)\big(q_a+(pB_{p-1})_1\big)\mpd$$
\end{Proposition}

Our next and last part is concerned with proving Theorem $0$.

\subsection{Convolutions of divided Bernoulli numbers and products of those with harmonic numbers.}

The work of \cite{LEV4} had allowed to find an expression for $A_{p-1}=(p-1)!$ modulo $p^4$. The formula expressed in terms of Bernoulli numbers was satisfactory in that it involved only a constant number of terms, independently from the prime $p$. This expression had also allowed to write the residue of a cubic convolution of order $(p-1)$ like follows.

\begin{equation*}
\mbt(p-1)=\Big\lbrace -6A_{p-1}+3p^3(\mc\mb(p-1))_1+3p^2(\mc\mb(p-1))_0-6p\mb_{p-1}-\frac{15}{4}p^3\mb_{p-3}\Big\rbrace_3\end{equation*}
\begin{equation*}
\qquad\qquad\qquad\qquad\qquad\qquad\qquad\qquad\qquad\qquad\qquad\qquad\qquad\qquad\qquad\mpu
\end{equation*}

Contrary to the generic case for $\mbt(p-1-2n)$ studied in $\S\,2.1$ and $\S\,2.3$, Gessel's identity does not provide any information on $\mbt(p-1)$. However interestingly, combined with Miki's identity it fully determines a convolution of divided Bernoulli numbers with a product of divided Bernoulli numbers and harmonic numbers. Indeed, setting $n=p-1$ in the Gessel identity (Theorem $0.4$ of $\S\,1.3$) yields (after using Wolstenholme's theorem imposing $\mh_{p-1}=0\mpd$):
\begin{equation}
\mbt(p-1)=\sum_{i=2}^{p-5}\mb_i\,b\mc\mb(p-1-j)+6\mb_{p-1}\as_{2,p-1}-\frac{9}{4}\mb_{p-3}\mpu
\end{equation}
When applying Miki's identity, the terms $\mbt(p-1)$ cancel each other, leaving only one convolution.
Moreover, a result by Jianqiang Zhao from \cite{ZO} claims that
$$\as_{2,p-1}=\frac{\mh_{p-1}}{p}\mpt$$
In \cite{LEV3}, combining our work and Zhao's result, we showed that
$$\as_{2,p-1}=-p(\mb_{2p-4}-2\mb_{p-3})\mpt$$
Plugging back into $(21)$ with $p\mb_{p-1}=1\mpu$ yields Theorem $0$ of $\S\,2.1$. It is interesting to note that by an application of Wolstenholme's theorem we have $\mh_{p-1-i}=\mh_i\;\mpu$. Then, by considering twice the convolutions of Corollary $0$, we get by a usual trick used many times in \cite{LEV4} the result of the corollary. \\\\

\hspace{0.5cm}\textsc{Email address:} \textit{clairelevaillant@yahoo.fr}

\end{document}